\documentclass[a4paper]{svproc}

\usepackage{amsmath, amssymb, graphicx, subcaption}
\usepackage{url}

\begin{document}
\setlength{\arraycolsep}{4.4pt}

\mainmatter              % start of a contribution
\title{Frequency-Domain Analysis of Controlled Euler--Bernoulli and Timoshenko Beams \\ with Attached Masses}
\titlerunning{Frequency-Domain Modeling of Controlled Flexible Beams}  % abbreviated title (for running head)
\author{Alexander Zuyev\inst{1,2} \and Julia Kalosha\inst{2}}
\authorrunning{Alexander Zuyev, Julia Kalosha} % abbreviated author list (for running head)
\tocauthor{Alexander Zuyev, Julia Kalosha}
\institute{Max Planck Institute for Dynamics of Complex Technical Systems\\ Sandtorstra{\ss}e 1, 39106 Magdeburg, Germany\\
\email{zuyev@mpi-magdeburg.mpg.de}
\and
Inst. of Applied Mathematics \& Mechanics, National Acad. of Sciences of Ukraine\\
\email{julykucher@gmail.com}}

\maketitle              % typeset the title of the contribution

\begin{abstract}
    This work focuses on the frequency-domain modeling of a control system with a flexible beam and a rigid body.
    A simply supported beam is equipped with a spring-loaded control actuator and possesses local damping effect.
    Using Hamilton's variational principle, the equations of motion are derived in the state space form taking into account interface conditions involving lumped control and local damping.
    The transfer functions are obtained for the Timoshenko and Euler--Bernoulli beam models with the output measurements provided by a point sensor.
    Comparative Bode plots are presented for the two beam models with different choices of output signals and damping coefficients.
\end{abstract}

\section{Introduction}

Modeling of mechanical systems involves balancing between complexity and accuracy.
In the theory of flexible beam vibrations with small transverse displacements, two widely accepted modeling approaches are distinguished: the Euler--Bernoulli and Timoshenko beam theories.
The former is suitable for modeling slender beams with sufficiently large aspect ratio (ratio of length to lateral dimensions).
When the effects of rotary inertia and shear deformation are negligible, the Euler–Bernoulli equation is often preferred due to its relative simplicity.
However, for shorter and stiffer beams with lower aspect ratios, shear deformation becomes significant~\cite{AR2021}. In such cases, the cross-section of the beam is not assumed to remain perpendicular to the neutral axis after the deformation, though it remains planar.
In this scenario, it is conventional to use a mathematical model in the form of the Timoshenko equations~\cite{E2019,T1921,T1974,KS2000}.

Comparative studies of different beam models have addressed various aspects, including spectral analysis~\cite{LBMF1996,TNC1953}, stability and control~\cite{AWS1997}, dynamic response~\cite{K2021,ZTS2020}, and analysis of mechanical characteristics~\cite{P2005,RB2007,W1995}.

The present work is devoted to modeling the motion of a damped flexible beam with a lumped rigid attachment.
We focus on the frequency domain representation and compare the numerical results obtained for the Timoshenko and Euler--Bernoulli beam models.
In Section~\ref{sec:TimBeamModel}, we derive a mathematical model in the form of the Timoshenko equations using Hamilton's variational principle.
Similar flexible structures within the Euler--Bernoulli beam theory, with rigid body attachment, were considered in our previous works~\cite{KZ2022,KZB2021,ZS05}.
An application of the variational principle to the derivation of the Timoshenko beam equations can be found, e.g., in~\cite{ZS2007}.
In Section~\ref{sec:TF.TimBeam}, we construct the transfer function thus presenting the damped Timoshenko model in the frequency domain.
This result opens up a wide range of possibilities for further investigation, in particular, analysis via frequency response techniques~\cite{B2013,JB2006}, stability analysis in the framework of frequency domain~\cite{GLO2017,M2024,WF2015}, and parameter identification~\cite{LHW1999}. % and reduced-order modeling~\cite{BBF2014,DO2014,KGA2021,LHW1999}.
Methods for constructing transfer functions for infinite-dimensional systems can be found in~\cite{CM2009,LA1993}.
Section~\ref{sec:EBmodel} presents the transfer function for the Euler--Bernoulli beam model describing the motion of an analogous mechanical system.
This allows us to compare the numerical results obtained for both analytical models, as discussed in Section~\ref{sec:numeric}.
An algorithm for deriving the transfer function for the Euler--Bernoulli beam was proposed in~\cite{ZG2023}, 
where a model with structural damping was investigated.
In the present paper, the dissipation effect is taken into account in the form of local viscous damping.

The principal features of the models considered in this paper are the presence of a point mass and local damping.
These peculiarities result in specific interface conditions that lead to additional free parameters and, consequently, require extra computations for the construction of the transfer function.

\section{Analytic Modeling of a Timoshenko Beam with a Rigid Body}\label{sec:TimBeamModel}

The mechanical system under consideration consists of a flexible beam of length $\ell$, simply supported at both ends.
The beam is assumed to be prismatic, homogeneous, and isotropic.
Directing the $x$-axis along the centerline of the beam, we introduce the function $w(x,t)$, which describes the transversal displacement of the point with coordinate $x\in[0,\ell]$ at time $t\geq0$,
and the function $\psi(x,t)$, which represents the shear angle of the cross section with respect to the main axis.
These functions are assumed to be of class $C^0$ w.r.t. $x$ and of class $C^2$ w.r.t. $t$.
Throughout this paper, we will use the following notations for the mechanical parameters of the system:
$\rho_0$ is the density (the mass per unit volume of the beam), $A$ is the cross section area, so $\rho=\rho_0A$ defines the linear density (the mass per unit length),
$I$ denotes the cross-sectional moment of inertia, so $I_\rho=\rho_0I$ is the mass moment of inertia of the cross section;
$E$ is Young's modulus; $G$ is the shear modulus (modulus of rigidity);
$K=kGA$, where $k$ is the shear correction factor that depends on the geometry of the beam's cross-section. We assume a rectangular cross-section, so $k=5/6$, see~\cite{DS2013,GT1991}.
Suppose that a rigid body (shaker) of mass $m$ is attached to the beam at a point $\ell_0 \in (0, \ell)$ by a hinge with spring stiffness $\varkappa$, which exerts a point load in the form of a force $F_0$.

We denote the spatial derivative by prime and the time derivative by dot.
The kinetic energy of the beam is presented as
$T = \frac12 \int\limits_0^\ell \big( \rho \dot w^2 + I_\rho \dot\psi^2 \big){\mathrm d}x$,
and the potential energy is
$U = \frac12 \int\limits_0^\ell \big( EI (\psi')^2 + K (w'-\psi)^2 \big){\mathrm d}x$.
The kinetic energy of the rigid body (shaker) is expressed as $\frac12 m \big(\dot w(\ell_0,t)\big)^2$,
the potential energy is $\frac12 \varkappa w^2(\ell_0,t)$, and the control force implemented by the actuator at $x=\ell_0$ is denoted by $F_0$.

Let $0\leq t_1<t_2$.
In order to derive mathematical model of the motion, we apply Hamilton's variational principle, which states
$\delta \int\limits_{t_1}^{t_2}(T-U+W_E){\mathrm d}t=0$ for any admissible variations
$\delta w(x,t)$, $\delta\psi(x,t)$, such that
$\delta w(0,t) = \delta w(\ell,t) = 0$, $\forall t\in[t_1,t_2]$,\;
$\delta w(x,t_1) = \delta w(x,t_2) = 0$,\; $\delta\psi(x,t_1) = \delta\psi(x,t_2) = 0$, $\forall x\in[0,\ell]$.
Here, $W_E$ describes the work of the control force $F_0$.
The admissible variations are assumed to be of class $C^1$.
Since the beam is simply supported, the boundary conditions
$w(0,t) = w(\ell,t) = 0$, $\forall t\in[t_1,t_2]$ are imposed.

After calculating the variation $\delta \int\limits_{t_1}^{t_2}(T-U+W_E){\mathrm d}t$ and applying the fundamental lemma of the calculus of variations, we obtain the following equations of motion:
\begin{equation}\label{eq:TB.eqs}
    \begin{array}{cc}
        \begin{aligned}
            \rho \ddot w(x,t) - K (w''(x,t) - \psi'(x,t)) & = 0, \\
            I_\rho \ddot\psi(x,t) - EI \psi''(x,t) - K (w'(x,t) - \psi(x,t)) & = 0,
        \end{aligned} & \quad
        x\neq\ell_0,
    \end{array}
\end{equation}
$w \in C^2\big(([0,\ell_0)\cup(\ell_0,\ell])\times[t_1,t_2]\big)$, $\psi(x,\cdot) \in C^2[t_1,t_2]$, $\forall x\in[0,\ell]$, $\psi(\cdot,t) \in C^1[0,\ell]$, $\forall t \in [t_1,t_2]$,
with the boundary conditions
\begin{equation}\label{eq:TB.BC}
    w(0,t) = w(\ell,t) = 0, \qquad \psi'(0,t) = \psi'(\ell,t) = 0
\end{equation}
and the interface conditions
\begin{equation}\label{eq:TB.IC}
    \begin{aligned}
        & w(\ell_0\!-\!0,t) = w(\ell_0\!+\!0,t), \;
        \psi(\ell_0\!-\!0,t) = \psi(\ell_0\!+\!0,t), \;
        \psi'(\ell_0\!-\!0,t) = \psi'(\ell_0\!+\!0,t), \\
        & F_0 - m \ddot w(\ell_0,t) - \varkappa w(\ell_0,t) - d \dot w(\ell_0,t) - K \Big( w'(\ell_0-0,t) - w'(\ell_0+0,t) \Big) = 0.
    \end{aligned}
\end{equation}
In~\eqref{eq:TB.IC}, we take into account local damping effect at a spring-mass attachment at the point $\ell_0$. The viscous damping is proportional to the transversal velocity of motion, with a coefficient $d>0$.

Assume that the system allows to measure output signals generated by a sensor attached at some point $\ell_k\in[0,\ell]$.
In this work, we consider two types of output signals:
the transversal displacement $y_1(t) = w(\ell_k,t)$ of the point $\ell_k\in[0,\ell]$ of the beam,
and the curvature of the beam axis at $\ell_k$ which can be defined as $\tilde y_2(t) = w''(\ell_k,t)$ for the Euler--Bernoulli beam.
As for the Timoshenko beam, at least two approaches for modeling the curvature can be found in literature: so called ``flexural'' (or ``mechanical''), and ``geometric'' curvature.
For detailed discussion of these models, we refer to~\cite{LCR2017} and references therein. In the present paper we consider the output of system~\eqref{eq:TB.eqs}--\eqref{eq:TB.IC} in the form of ``mechanical'' curvature which is modeled as $y_2(t)=\psi'(\ell_k,t)$.

\section{Transfer Function for the Damped Timoshenko Beam Model}\label{sec:TF.TimBeam}
In this section, we propose an algorithm to obtain the frequency-domain representation for the dynamical system derived above.
We assume zero initial conditions and apply the Laplace transform w.r.t. $t$ to input, output and unknown functions in system \eqref{eq:TB.eqs}--\eqref{eq:TB.IC}.
The Laplace transform operator is denoted by ${\cal L}: g(x,t)\mapsto G(x,s)=\int\limits_0^\infty g(x,t) e^{-st}{\mathrm d}t$, where $s\in{\mathbb C}$ is a parameter.
Denote $U(s)={\cal L}(F_0(t))$, $Y(s)={\cal L}(y(t))$, $W(x,s)={\cal L}(w(x,t))$, $\Psi(x,s)={\cal L}(\psi(x,t))$,
%then ${\cal L}(\dot w(x,t))=sW(x,s)$, ${\cal L}(\ddot w(x,t))=s^2W(x,s)$, ${\cal L}(\dot\psi(x,t))=s\Psi(x,s)$, ${\cal L}(\ddot \psi(x,t))=s^2\Psi(x,s)$.
%Denoting
%$U(s) = \int\limits_0^\infty F_0(t) e^{-st}{\mathrm d}t$,\;
%$Y(s) = \int\limits_0^\infty y(t) e^{-st}{\mathrm d}t$,\;
%$W(x,s) = \int\limits_0^\infty w(x,t) e^{-st}{\mathrm d}t$,\;
%$\Psi(x,s) = \int\limits_0^\infty \psi(x,t) e^{-st}{\mathrm d}t$ with parameter $s\in{\mathbb C}$,
then we obtain the following system:
\begin{equation}\label{eq:TB.LT.sys}
    \begin{array}{cc}
        \begin{aligned}
            K \big(W''(x,s) - \Psi'(x,s) \big) - \rho s^2\, W(x,s) & = 0, \\
            EI\,\Psi''(x,s) + K \big(W'(x,s) - \Psi(x,s) \big) - I_\rho s^2\, \Psi(x,s) & = 0,
            \end{aligned} & \quad
        x\neq\ell_0,
    \end{array}
\end{equation}
\begin{equation}\label{eq:TB.LT.BC}
    W(0,s) = W(\ell,s) = 0, \quad
    \Psi'(0,s) = \Psi'(\ell,s) = 0,
\end{equation}
\begin{equation}\label{eq:TB.LT.IC}
    \begin{aligned}
        & W(\ell_0-0,s) = W(\ell_0+0,s), \\
        & \Psi(\ell_0-0,s) = \Psi(\ell_0+0,s), \quad \Psi'(\ell_0-0,s) = \Psi'(\ell_0+0,s), \\
        & K\,\big(W'(\ell_0-0,s) - W'(\ell_0+0,s)\big) + (ms^2 + \varkappa + d s)\, W(\ell_0,s) = U(s).
    \end{aligned}
\end{equation}
A series of transformations applied to system~\eqref{eq:TB.LT.sys} leads to the forth-order ODE with parameter $s$:
\begin{equation}\label{eq:TB.ODE}
    KEI\,W^{(4)}(x,s) - \big( \rho EI+KI_\rho \big)s^2\, W''(x,s) + \big(K+I_\rho\, s^2\big)\rho s^2\, W(x,s) = 0.
\end{equation}
Here and throughout the text, the upper index in parentheses denotes the spatial derivative of the corresponding order, i.e., $W^{(j)}(x,s)=\dfrac{{\mathrm d}^j}{{\mathrm d}x^j}\,W(x,s)$, $j=\overline{1,4}$.

In what follows throughout this section, we assume the following:
\begin{equation}\label{eq:sI1}
s\not\in{\cal I} = \left\{0,\; \pm i\sqrt{\dfrac{K}{I_\rho}},\; \pm\frac{2K\sqrt{\rho EI}}{\rho EI-KI_\rho} \right\}.
\end{equation}

Denoting
\begin{multline}\label{eq:Wvec}
    W_j=W^{(j)}(x,s), \; W_j^0=W^{(j)}(0,s), \; W_j^\ell=W^{(j)}(\ell,s), \; j=\overline{0,3}, \\
    \bar W = \begin{pmatrix} W_0, & W_1, & W_2, & W_3 \end{pmatrix}^\top, \; \bar W^0 = \bar W(0,s)$, $\bar W^\ell = \bar W(\ell,s),
\end{multline}
we represent equation~\eqref{eq:TB.ODE} as the following first-order system:
$
    \frac{\mathrm d}{{\mathrm d}x}\,\bar W = {\cal A}\, \bar W,
$
$${\cal A}(s) = \left(\begin{array}{cccc}
        0  & 1 & 0  & 0 \\
        0  & 0 & 1  & 0 \\
        0  & 0 & 0  & 1 \\
        -b & 0 & 2a & 0
    \end{array}\right), \quad
a(s) = \frac{s^2}2\left(\frac\rho K + \frac{I_\rho}{EI}\right), \quad
b(s) = \frac{(K+I_\rho s^2)\rho s^2}{KEI},$$
with the solution expressed through the matrix exponential in the following form:
\begin{equation}\label{eq:ODEsol}
    \bar W(x,s) = \left\{\begin{array}{ll}
        e^{x{\cal A}}\, \bar W^0,           & x\in[0,\ell_0], \\
        e^{(x-\ell){\cal A}}\, \bar W^\ell, & x\in(\ell_0,\ell],
    \end{array}\right.
\end{equation}
where
$e^{x{\cal A}}(s) = \dfrac1{\lambda_1^2-\lambda_2^2}
    \left(\begin{array}{cccc}
        z_1                          & z_2                          & z_3 & z_4 \\
        -\lambda_1^2\lambda_2^2\,z_4 & z_1                          & z_5 & z_3 \\
        -\lambda_1^2\lambda_2^2\,z_3 & -\lambda_1^2\lambda_2^2\,z_4 & z_6 & z_5 \\
        -\lambda_1^2\lambda_2^2\,z_5 & -\lambda_1^2\lambda_2^2\,z_3 & z_7 & z_6
    \end{array}\right)$, \\
$\lambda_1(s) = \sqrt{a+\sqrt{a^2-b}}$, $\lambda_2(s) = \sqrt{a-\sqrt{a^2-b}}$,
\begin{equation}\label{eq:TB.zfu}
    \begin{aligned}
        & z_1(x,s) = \lambda_1^2\cosh\lambda_2x-\lambda_2^2\cosh\lambda_1x, \\
        & z_2(x,s) = \frac1{\lambda_1\lambda_2}(\lambda_1^3\sinh\lambda_2x-\lambda_2^3\sinh\lambda_1x), \;
            z_3(x,s) = \cosh\lambda_1x-\cosh\lambda_2x, \\
        & z_4(x,s) = \frac1{\lambda_1\lambda_2}(\lambda_2\sinh\lambda_1x-\lambda_1\sinh\lambda_2x), \;
            z_5(x,s) = \lambda_1\sinh\lambda_1x-\lambda_2\sinh\lambda_2x, \\
        & z_6(x,s) = \lambda_1^2\cosh\lambda_1x-\lambda_2^2\cosh\lambda_2x, \;
            z_7(x,s) = \lambda_1^3\sinh\lambda_1x-\lambda_2^3\sinh\lambda_2x.
    \end{aligned}
\end{equation}

Observe that the boundary conditions~\eqref{eq:TB.LT.BC} immediately imply that
$W_0^0=W_2^0=0$,\; $W_0^\ell=W_2^\ell=0$, so the solution of~\eqref{eq:TB.LT.sys}--\eqref{eq:TB.LT.BC} reads as
\begin{equation}\label{eq:BVPsolW}
    (\lambda_1^2-\lambda_2^2)\,W(x,s) =
    \left\{\begin{array}{ll}
        z_2(x,s)\,W_1^0 + z_4(x,s)\,W_3^0,                 & \;\; x\in[0,\ell_0], \\
        z_2(x-\ell,s)\,W_1^\ell + z_4(x-\ell,s)\,W_3^\ell, & \;\; x\in(\ell_0,\ell],
    \end{array}\right.
\end{equation}
\begin{equation}\label{BVPsolP}
    (\lambda_1^2-\lambda_2^2)\,\Psi(x,s) =
    \left\{\begin{array}{ll}       
        \begin{aligned}
            u\Big(vz_1(x,s)-\lambda_1^2\lambda_2^2\,z_3(x,s)\Big)W_1^0 \\[.25pt]
            + \Big(vz_3(x,s)+z_6(x,s)\Big)W_3^0,
        \end{aligned} & \;\; x\in[0,\ell_0], \\
        \begin{aligned}
            u\Big(vz_1(x-\ell,s)-\lambda_1^2\lambda_2^2\,z_3(x-\ell,s)\Big)W_1^\ell \\[.25pt]
            + \Big(vz_3(x-\ell,s)+z_6(x-\ell,s)\Big)W_3^\ell,
        \end{aligned} & \;\; x\in(\ell_0,\ell],
    \end{array}\right.
\end{equation}
where $u(s) = \frac{EI}{K+I_\rho s^2}$,\; $v(s) = \frac K{EI}-\frac{\rho s^2}K$.

In order to eliminate the remaining constants of integration we exploit the interface conditions~\eqref{eq:TB.LT.IC} that result in the algebraic system \linebreak
$M \begin{pmatrix} W_1^0, & W_3^0, & W_1^\ell, & W_3^\ell \end{pmatrix}^\top = \begin{pmatrix} 0, & 0, & 0, & \frac{\lambda_1^2-\lambda_2^2}K\, U(s) \end{pmatrix}^\top$
with the matrix
\begin{equation}\label{eq:TB.mM}
    M(s) = \left(\begin{array}{cccc}
        z_2(\ell_0,s) & z_4(\ell_0,s) & -z_2(\ell_0-\ell,s) & -z_4(\ell_0-\ell,s) \\
        M_{21}        & M_{22}        & M_{23}              & M_{24} \\
        M_{31}        & M_{32}        & M_{33}              & M_{34} \\
        M_{41}        & M_{42}        & -z_1(\ell_0-\ell,s) & -z_3(\ell_0-\ell,s)
    \end{array}\right),
\end{equation}
where
$M_{21} = vz_1(\ell_0,s)-\lambda_1^2\lambda_2^2\,z_3(\ell_0,s)$, \;
$M_{22} = vz_3(\ell_0,s)+z_6(\ell_0,s)$, \\
$M_{23} = \lambda_1^2\lambda_2^2\,z_3(\ell_0-\ell,s)-vz_1(\ell_0-\ell,s)$, \;
$M_{24} = -z_6(\ell_0-\ell,s)-vz_3(\ell_0-\ell,s)$, \\
$M_{31} = v_1z_2(\ell_0,s)+\lambda_1^2\lambda_2^2\,z_4(\ell_0,s)$, \;
$M_{32} = v_1z_4(\ell_0,s)-z_5(\ell_0,s)$, \\
$M_{33} = -\lambda_1^2\lambda_2^2\,z_4(\ell_0-\ell,s)-v_1z_2(\ell_0-\ell,s)$, \;
$M_{34} = z_5(\ell_0-\ell,s)-v_1z_4(\ell_0-\ell,s)$, \\
$M_{41} = z_1(\ell_0,s)+v_2z_2(\ell_0,s)$, \;
$M_{42} = z_3(\ell_0,s)+v_2z_4(\ell_0,s)$, \\
%$$M_{21} = vz_1(\ell_0,s)-\lambda_1^2\lambda_2^2\,z_3(\ell_0,s), \quad M_{31} = v_1z_2(\ell_0,s)+\lambda_1^2\lambda_2^2\,z_4(\ell_0,s),$$
%$$M_{22} = vz_3(\ell_0,s)+z_6(\ell_0,s), \quad M_{32} = v_1z_4(\ell_0,s)-z_5(\ell_0,s),$$
%$$M_{23} = \lambda_1^2\lambda_2^2\,z_3(\ell_0-\ell,s)-vz_1(\ell_0-\ell,s), \; M_{33} = -\lambda_1^2\lambda_2^2\,z_4(\ell_0-\ell,s)-v_1z_2(\ell_0-\ell,s),$$
%$$M_{24} = -z_6(\ell_0-\ell,s)-vz_3(\ell_0-\ell,s), \quad M_{34} = z_5(\ell_0-\ell,s)-v_1z_4(\ell_0-\ell,s),$$
%$$M_{41} = z_1(\ell_0,s)+v_2z_2(\ell_0,s), \quad M_{42} = z_3(\ell_0,s)+v_2z_4(\ell_0,s),$$
$v_1(s) = \frac{\rho s^2}{K}$, \; $v_2(s) = \frac{ms^2+ds+\varkappa}{K}$.

If the Laplace parameter $s$ is not an eigenvalue of the differential operator corresponding to system~\eqref{eq:TB.eqs}--\eqref{eq:TB.IC}, and $s$ satisfies assumption~\eqref{eq:sI1}, then matrix $M$ is nonsingular and the unique solution of~\eqref{eq:TB.LT.sys}--\eqref{eq:TB.LT.IC} is obtained by substituting
$\begin{pmatrix} W_1^0, & W_3^0, & W_1^\ell, & W_3^\ell \end{pmatrix}^\top = \frac{\lambda_1^2-\lambda_2^2}K\, U \begin{pmatrix} M_{14}^{-1}, & M_{24}^{-1}, & M_{34}^{-1}, & M_{44}^{-1} \end{pmatrix}^\top$
%$W_1^0 = M_{14}^{-1} \frac{\lambda_1^2-\lambda_2^2}K\, U(s)$, $W_3^0 = M_{24}^{-1} \frac{\lambda_1^2-\lambda_2^2}K\, U(s)$, $W_1^\ell = M_{34}^{-1} \frac{\lambda_1^2-\lambda_2^2}K\, U(s)$ and $W_3^\ell = M_{44}^{-1} \frac{\lambda_1^2-\lambda_2^2}K\, U(s)$
into~\eqref{eq:BVPsolW},~\eqref{BVPsolP}.
%$$W(x,s) = \left\{\begin{array}{ll}
%        \Big(z_2(x,s) M_{14}^{-1} + z_4(x,s) M_{24}^{-1}\Big) \frac{U(s)}{K},           & \;\; x\in[0,\ell_0], \\
%        \Big(z_2(x-\ell,s) M_{34}^{-1} + z_4(x-\ell,s) M_{44}^{-1}\Big) \frac{U(s)}{K}, & \;\; x\in(\ell_0,\ell].
%    \end{array}\right.$$

The transfer function is calculated as $H(s) = \frac{Y(s)}{U(s)}$.
Given that the output signal is provided in the form $Y_1(s)=W(\ell_k,s)$, or $Y_2(s)=\Psi'(\ell_k,s)$, we obtain the frequency-domain representation as proposed below.
\begin{proposition}\label{le:TB.TF}
    Let $S=\{s\in{\mathbb C}\mid \det(M(s))\neq0\}$.
    The transfer functions $H_1(s)$ and $H_2(s)$ of the control system~\eqref{eq:TB.eqs}--\eqref{eq:TB.IC}
    with outputs $y_1(t) = w(\ell_k,t)$ and $y_2(t) = \psi'(\ell_k,t)$, respectively, are defined for all $s\in S\setminus{\cal I}$ as follows:
    $$H_1(s) = \left\{\begin{array}{ll}
        \frac1K \Big(z_2(\ell_k,s) M_{14}^{-1} + z_4(\ell_k,s) M_{24}^{-1}\Big),           & \ell_k\in[0,\ell_0], \\
        \frac1K \Big(z_2(\ell_k-\ell,s) M_{34}^{-1} + z_4(\ell_k-\ell,s) M_{44}^{-1}\Big), & \ell_k\in(\ell_0,\ell],
        \end{array}\right.$$
    $$H_2(s) = \left\{\begin{array}{ll}
        \begin{aligned}
            \frac1{K^2} \Big[\big(-\lambda_1^2\lambda_2^2\,K\,z_4(\ell_k,s) - \rho s^2\,z_2(\ell_k,s)\big) M_{14}^{-1} \Big. \\
            + \Big.\big(K\,z_5(\ell_k,s) - \rho s^2\,z_4(\ell_k,s)\big) M_{24}^{-1}\Big],
        \end{aligned}
        & \ell_k\in[0,\ell_0], \\
        \begin{aligned}
            \frac1{K^2} \Big[\big(-\lambda_1^2\lambda_2^2\,K\,z_4(\ell_k-\ell,s) - \rho s^2\,z_2(\ell_k-\ell,s)\big) M_{34}^{-1} \Big. \\
            + \Big.\big(K\,z_5(\ell_k-\ell,s) - \rho s^2\,z_4(\ell_k-\ell,s)\big) M_{44}^{-1}\Big],
        \end{aligned}
        & \ell_k\in(\ell_0,\ell]
        \end{array}\right.$$
    where the components $z_2(x,s)$, $z_4(x,s)$, $z_5(x,s)$ are specified in~\eqref{eq:TB.zfu}, and $M_{j4}^{-1}$, $j=\overline{1,4}$, denote the elements of the inverse matrix corresponding to~\eqref{eq:TB.mM}.
\end{proposition}

\section{Damped Euler--Bernoulli Beam Model in the Frequency Domain}\label{sec:EBmodel}

Consider the mechanical setup described in Section~\ref{sec:TimBeamModel}.
In this section, we derive its mathematical model in the form of the Euler--Bernoulli equation and obtain the corresponding transfer function.
Within the scope of Euler--Bernoulli linear theory, we employ the assumption that the rotation of the beam's element is insignificant compared to the transversal displacement, so the shear deformation of the cross-section is neglected.

The kinetic and potential energies of the beam are presented, respectively, as
$T = \int\limits_0^\ell \rho(\dot w(x,t))^2{\mathrm d}x + m(\dot w(\ell_0,t))^2$ and
$U = \int\limits_0^\ell  EI(w''(x,t))^2{\mathrm d}x + \varkappa(w^2(\ell_0,t))$.
Applying Hamilton's principle yields the equation of motion:
\begin{equation}\label{eq:EB.eq}
    \ddot w(x,t) + \frac{EI}\rho\, w^{(4)}(x,t) = 0, \quad x\neq\ell_0,
\end{equation}
subject to the boundary and interface conditions:
\begin{equation}\label{eq:EB.BC}
    w(0,t) = w(\ell,t) = 0, \quad w''(0,t) = w''(\ell,t) = 0,
\end{equation}
\begin{equation}\label{eq:EB.IC}
    \begin{aligned}
        & w^{(j)}(\ell_0-0,t) = w^{(j)}(\ell_0+0,t),\: j=\overline{0,2}, \\
        & \big(m\ddot w + \varkappa w + d\dot w\big)\Big|_{x=\ell_0} = EI\left(w'''\Big|_{x=\ell_0-0} - w'''\Big|_{x=\ell_0+0}\right) + F_0.
    \end{aligned}
\end{equation}

We assume that this system admits the same types of output signals as the Timoshenko model studied in the previous section.
The procedure of obtaining the transfer function for system~\eqref{eq:EB.eq}--\eqref{eq:EB.IC} is similar to the one described above.
Therefore, we omit the detailed description and instead outline the key steps.

After applying the Laplace transform w.r.t. $t$ to equations~\eqref{eq:EB.eq}--\eqref{eq:EB.IC} with vanishing initial conditions, we obtain the following system:
\begin{equation}\label{eq:EB.LT.sys}
    s^2\, W(x,s) + \frac{EI}\rho W^{(4)}(x,s) = 0, \quad x\neq\ell_0,
\end{equation}
\begin{equation}\label{eq:EB.LT.BC}
    W(0,s) = W(\ell,s) = 0, \quad W''(0,s) = W''(\ell,s) = 0,
\end{equation}
\begin{equation}\label{eq:EB.LT.IC}
    \begin{aligned}
        & W^{(j)}(\ell_0-0,s) = W^{(j)}(\ell_0+0), \; j=\overline{0,2}, \\
        & (ms^2+ds+\varkappa) W(\ell_0,s) = EI\big(W'''(\ell_0-0,s) - W'''(\ell_0+0,s)\big) + U(s),
    \end{aligned}
\end{equation}
where $s$ is the Laplace variable.

The parametric ODE~\eqref{eq:EB.LT.sys} has a solution of the form~\eqref{eq:ODEsol}, where we replace ${\cal A}$ with $\tilde{\cal A}$,
\begin{equation*}
    \tilde{\cal A}(s) = \left(\begin{array}{cccc}
        0                    & 1 & 0 & 0 \\
        0                    & 0 & 1 & 0 \\
        0                    & 0 & 0 & 1 \\
        -\frac{\rho s^2}{EI} & 0 & 0 & 0
    \end{array}\right), \;
    e^{x\tilde{\cal A}}(s) = \frac12 \left(\begin{array}{cccc}
        \tilde z_1           & \tilde z_2           & \tilde z_3           & \tilde z_4 \\
        \gamma^4\,\tilde z_4 & \tilde z_1           & \tilde z_2           & \tilde z_3 \\
        \gamma^4\,\tilde z_3 & \gamma^4\,\tilde z_4 & \tilde z_1           & \tilde z_2 \\
        \gamma^4\,\tilde z_2 & \gamma^4\,\tilde z_3 & \gamma^4\,\tilde z_4 & \tilde z_1
    \end{array}\right),
\end{equation*}
where $\gamma=\gamma(s)$ denotes the principal value of $\sqrt[4]{-\frac{\rho s^2}{EI}}$, and
\begin{equation}\label{eq:EB.zfu}
    \begin{array}{ll}
        \tilde z_1(x,s) = \cosh\gamma x + \cos\gamma x,                        & \tilde z_3(x,s) = \frac1{\gamma^2} \big(\cosh\gamma x - \cos\gamma x\big), \\
        \tilde z_2(x,s) = \frac1\gamma \big(\sinh\gamma x + \sin\gamma x\big), & \tilde z_4(x,s) = \frac1{\gamma^3} \big(\sinh\gamma x - \sin\gamma x\big).
    \end{array}
\end{equation}
Further, we continue using the notation~\eqref{eq:Wvec}.
Boundary conditions~\eqref{eq:EB.LT.BC} imply
$W_0^0 = W_0^\ell = 0$,\; $W_2^0 = W_2^\ell = 0$.
The algebraic system for eliminating the remaining coefficients, obtained from interface conditions~\eqref{eq:EB.LT.IC}, reads as \linebreak
$\tilde M \begin{pmatrix} W_1^0, & W_3^0, & W_1^\ell, & W_3^\ell \end{pmatrix}^\top = \begin{pmatrix} 0, & 0, & 0, & \frac{2\,U(s)}{EI} \end{pmatrix}^\top$,
where
\begin{equation}\label{eq:EB.mM}
    \tilde M(s) = \left(\begin{array}{cccc}
        \tilde z_2(\ell_0,s)         & \tilde z_4(\ell_0,s) & -\tilde z_2(\ell_0-\ell,s)         & -\tilde z_4(\ell_0-\ell,s) \\
        \tilde z_1(\ell_0,s)         & \tilde z_3(\ell_0,s) & -\tilde z_1(\ell_0-\ell,s)         & -\tilde z_3(\ell_0-\ell,s) \\
        \gamma^4\tilde z_4(\ell_0,s) & \tilde z_2(\ell_0,s) & -\gamma^4\tilde z_4(\ell_0-\ell,s) & -\tilde z_2(\ell_0-\ell,s) \\
        \tilde M_{41}                & \tilde M_{42}        & \gamma^4\tilde z_3(\ell_0-\ell,s)  & \tilde z_1(\ell_0-\ell,s)
    \end{array}\right),
\end{equation}
$\tilde M_{41} = \tilde v\,\tilde z_2(\ell_0,s)-\gamma^4\,\tilde z_3(\ell_0,s)$,\;
$\tilde M_{42} = \tilde v\,\tilde z_4(\ell_0,s)-\tilde z_1(\ell_0,s)$,\;
$\tilde v(s) = \frac{m s^2+ds+\varkappa}{EI}$.
The inverse matrix $\tilde M^{-1} = (\tilde M_{ij}^{-1})_{i,j=1}^4$ exists provided that the parameter $s$ is not an eigenvalue of the differential operator corresponding to system~\eqref{eq:EB.eq}--\eqref{eq:EB.IC}.

\begin{proposition}\label{le:EB.TF}
    The transfer functions $\tilde H_1$ and $\tilde H_2$ of the control system~\eqref{eq:EB.eq}--\eqref{eq:EB.IC},
    with outputs $y_1(t) = w(\ell_k,t)$ and $\tilde y_2(t) = w''(\ell_k,t)$, respectively, are defined for all $s\in{\mathbb C}$ such that $\det(\tilde M(s))\neq0$, as follows:
    $$\tilde H_1(s) = \left\{\begin{array}{ll}
        \frac1{EI} \Big(\tilde z_2(\ell_k,s) \tilde M_{14}^{-1} + \tilde z_4(\ell_k,s) \tilde M_{24}^{-1}\Big),           & \ell_k\in[0,\ell_0], \\
        \frac1{EI} \Big(\tilde z_2(\ell_k-\ell,s) \tilde M_{34}^{-1} + \tilde z_4(\ell_k-\ell,s) \tilde M_{44}^{-1}\Big), & \ell_k\in(\ell_0,\ell],
        \end{array}\right.$$
    $$\tilde H_2(s) = \left\{\begin{array}{ll}
        \frac1{EI} \Big(\gamma^4 \tilde z_4(\ell_k,s) \tilde M_{14}^{-1} + \tilde z_2(\ell_k,s) \tilde M_{24}^{-1}\Big),    & \ell_k\in[0,\ell_0], \\
        \frac1{EI} \Big(\gamma^4 \tilde z_4(\ell_k-\ell,s) M_{34}^{-1} + \tilde z_2(\ell_k-\ell,s) \tilde M_{44}^{-1}\Big), & \ell_k\in(\ell_0,\ell],
        \end{array}\right.$$
    where the components $\tilde z_2(x,s)$, $\tilde z_4(x,s)$ are determined in~\eqref{eq:EB.zfu}, and $\tilde M_{j4}^{-1}$, $j=\overline{1,4}$,
    are the elements of the matrix inverse to~\eqref{eq:EB.mM}.
\end{proposition}

\section{Comparative Analysis of the Euler--Bernoulli and Timoshenko Beam Models in the Frequency Domain}\label{sec:numeric}

In this section, we present comparative Bode plots for the two oscillating system models investigated above.
For the numeric simulations, we use the following values of mechanical parameters:
\begin{align*}
    & \ell=1.905~{\rm m}, \; \ell_0=1.4~{\rm m}, \; A=2.25~{\rm cm^2}, \; m=0.045~{\rm kg}, \; \varkappa=2630~{\rm N/m}, \\
    & \rho_0=2700~{\rm kg/m^3}, \; E=69~{\rm GPa}, \; G=25.5~{\rm GPa}, \; I=1.6875\cdot10^{-10}~{\rm m^4}.
\end{align*}

\begin{figure}[!ht]\centering
    % ===== LEFT COLUMN =====
    \begin{minipage}[b]{0.47\linewidth} \centering
        \begin{subfigure}{\linewidth} \centering
            \includegraphics[width=\linewidth]{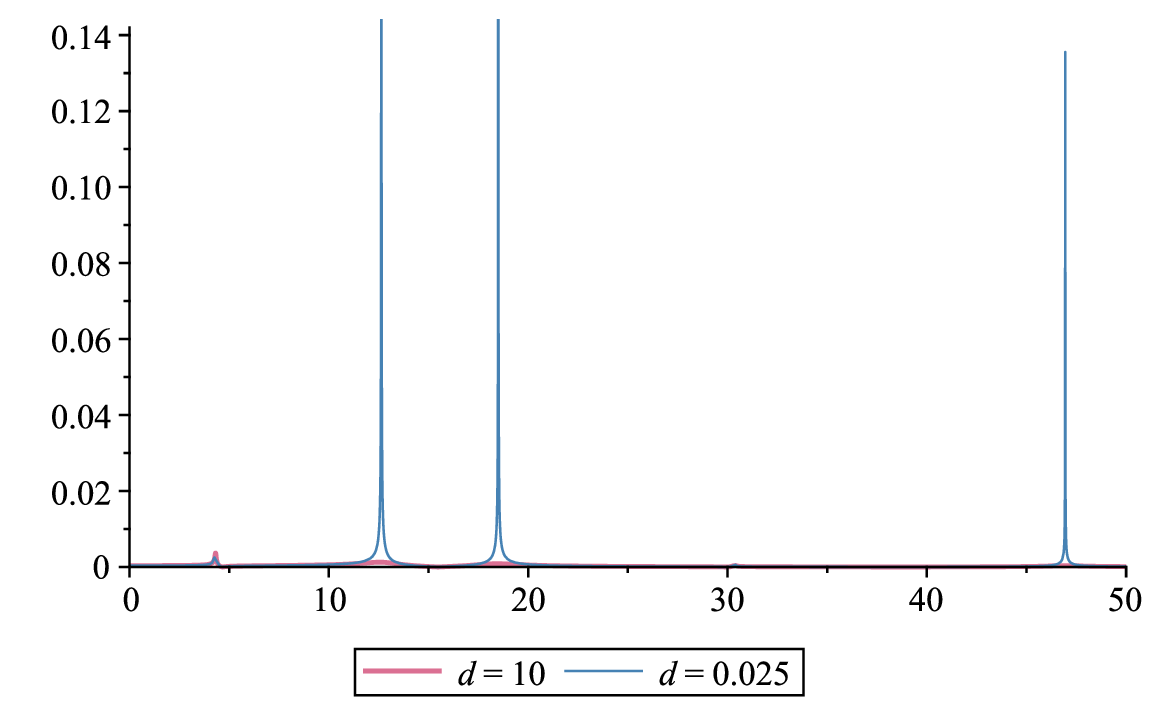}
            \caption{$|H_1(2i\pi\nu)|$, Timoshenko model.}\label{fig:TB.op1}
        \end{subfigure}
        \begin{subfigure}{\linewidth} \centering
            \includegraphics[width=\linewidth]{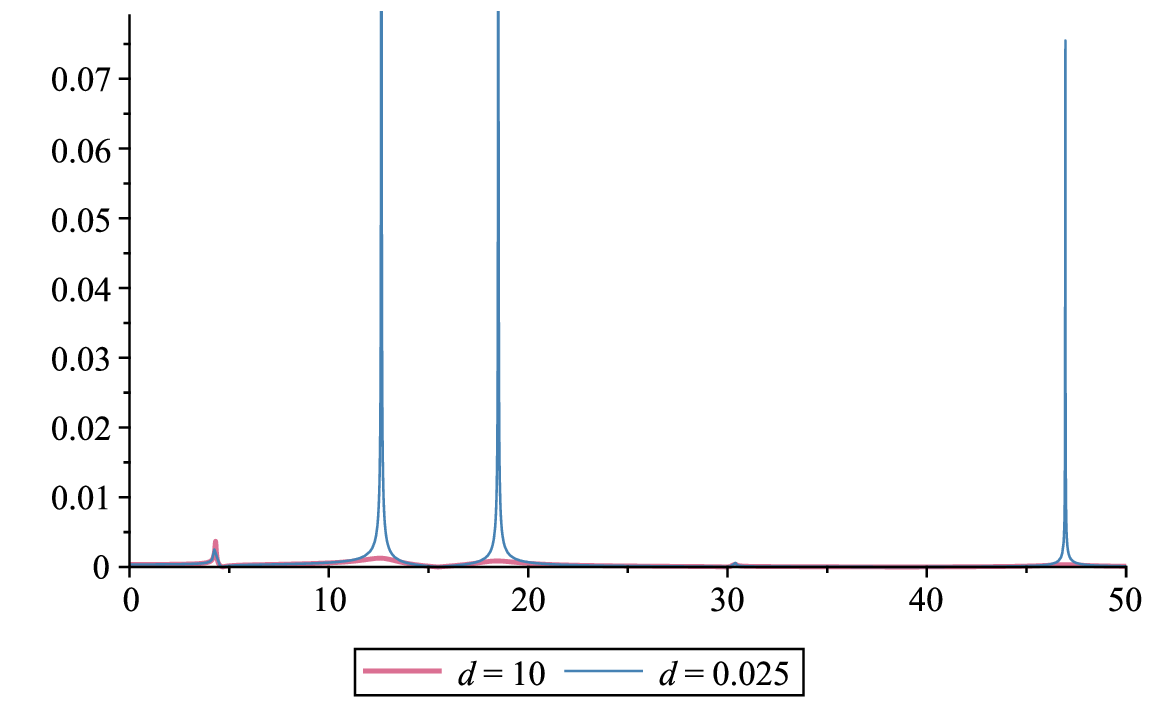}
            \caption{$|\tilde H_1(2i\pi\nu)|$, Euler--Bernoulli model.}\label{fig:EB.op1}
        \end{subfigure}
        \caption{Magnitude Bode plots for systems with different damping coefficients. Output $y_1(t)$. \phantom{Lorem ipsum dolor}}\label{fig:op1}
    \end{minipage} \hfill
    % ===== RIGHT COLUMN =====
    \begin{minipage}[b]{0.47\linewidth} \centering
        \begin{subfigure}{\linewidth} \centering
            \includegraphics[width=\linewidth]{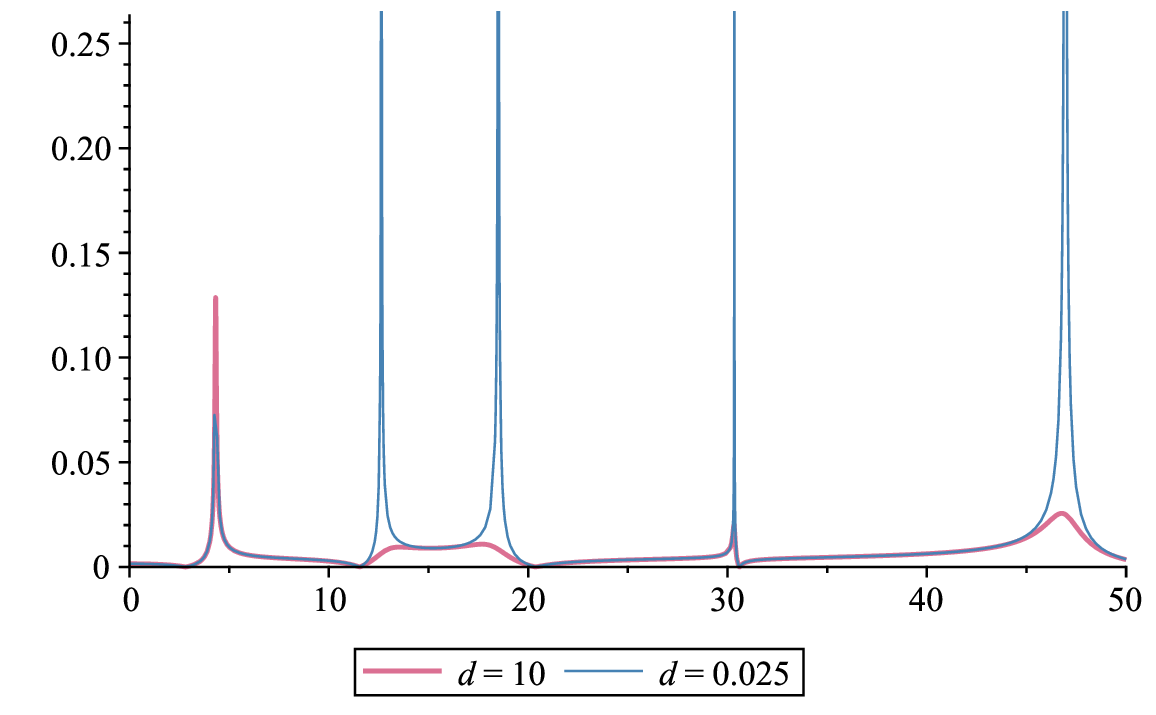}
            \caption{$|H_2(2i\pi\nu)|$, Timoshenko model.}\label{fig:TB.op2}
        \end{subfigure}
        \begin{subfigure}{\linewidth} \centering
            \includegraphics[width=\linewidth]{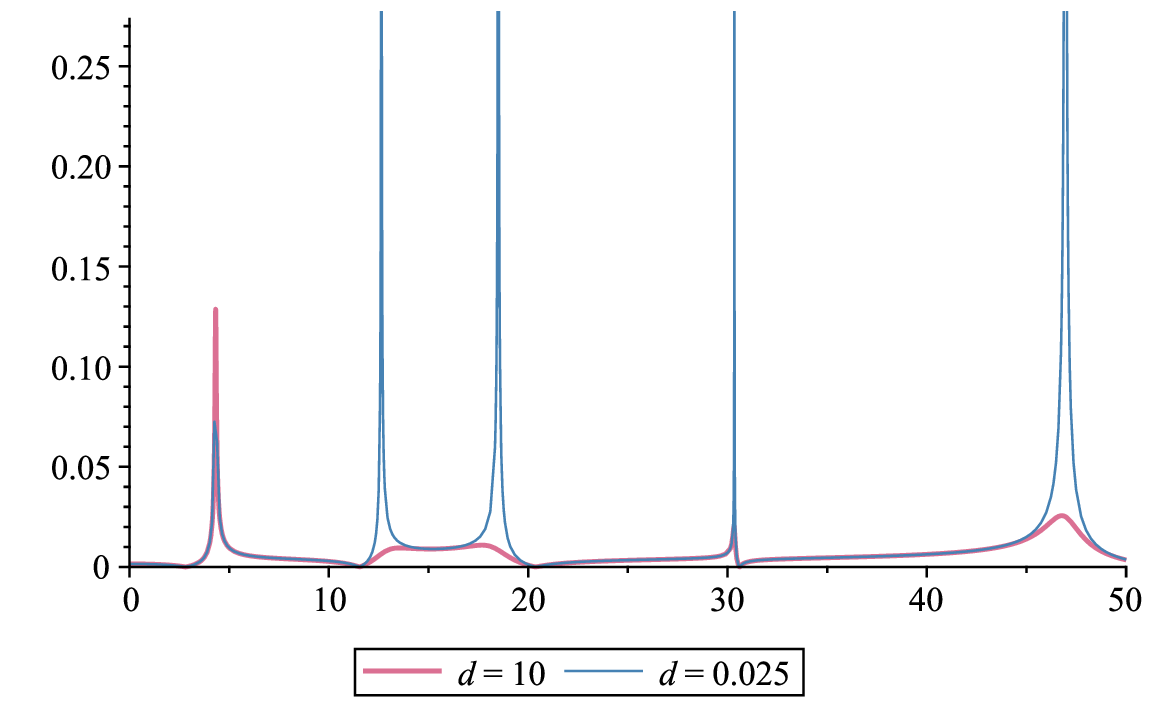}
            \caption{$|\tilde H_2(2i\pi\nu)|$, Euler--Bernoulli model.}\label{fig:EB.op2}
        \end{subfigure}
        \caption{Magnitude Bode plots for systems with different damping coefficients. Outputs $y_2(t)$~(top) and $\tilde y_2(t)$~(bottom).}\label{fig:op2}
    \end{minipage} \vspace{-1ex}
\end{figure}

The following figures present magnitude Bode plots for the Timoshenko and Euler--Bernoulli beam models with rigid attachments.
The plots are obtained as magnitudes of transfer functions derived in Propositions~\ref{le:TB.TF} and~\ref{le:EB.TF} calculated for argument values $2i\pi\nu$ with $\nu\in[0,50]$ standing for the frequency in~Hz.
Namely, Fig.~\ref{fig:op1} is obtained for the Timoshenko beam model (Fig.~\ref{fig:TB.op1}) and the Euler--Bernoulli model (Fig.~\ref{fig:EB.op1}) in the case of the displacement output signal $y_1(t)=w(\ell_k,t)$,
while Fig.~\ref{fig:op2} is for curvature-type output signals: $y_2(t)=\psi'(\ell_k,t)$ for the Timoshenko model on Fig.~\ref{fig:TB.op2} and $\tilde y_2(t)=w''(\ell_k,t)$ for the Euler--Bernoulli model on Fig.~\ref{fig:EB.op2}. All output signals are evaluated at the point $\ell_k=\ell_0$. 
These plots illustrate the effect of variations in the damping coefficient $d$ on the transfer functions magnitudes for both flexible beam models, under the two types of output described above, in low-frequency range $\nu\in[0,50]$.

In Table~\ref{tab:freq}, we present the modal frequencies as they depend on the damping coefficient $d$.
For that, we compute the local maxima points $\nu_j$ and $\tilde\nu_j$ on the real line for the magnitudes of the transfer functions $|H_{1,2}(2i\pi\nu)|$~(Timoshenko beam model) and $|\tilde H_{1,2}(2i\pi\nu)|$~(Euler--Bernoulli beam), respectively.
Table~\ref{tab:freq.op1} contains results obtained for both systems with displacement outputs $y_1=w(\ell_0,t)$, while
Table~\ref{tab:freq.op2} presents the corresponding frequencies for systems with curvature-type outputs: $y_2=\psi'(\ell_0,t)$ for the Timoshenko model and $\tilde y_2=w''(\ell_0,t)$ for the Euler--Bernoulli model.

\begin{table}[!ht]
    \caption{Modal frequencies}\label{tab:freq}
    \vspace{-1ex}
    \begin{subtable}{0.5\linewidth}
        \caption{Output $y_1(t)$.}\label{tab:freq.op1} \centering
        \begin{tabular}{|c||c|c|c|} \hline
                          & $d=0.025$ & $d=1$     & $d=10$    \\ \hline
            $\nu_1$       & $4.3185$  & $4.3185$  & $4.3182$  \\ 
            $\tilde\nu_1$ & $4.3186$  & $4.3186$  & $4.3183$  \\ \hline
            $\nu_2$       & $12.6358$ & $12.6352$ & $12.5700$ \\ 
            $\tilde\nu_2$ & $12.6360$ & $12.6354$ & $12.5705$ \\ \hline
            $\nu_3$       & $18.5019$ & $18.5012$ & $18.4411$ \\ 
            $\tilde\nu_3$ & $18.5025$ & $18.5019$ & $18.4415$ \\ \hline
            $\nu_4$       & $30.3447$ & $30.3447$ & $30.3444$ \\ 
            $\tilde\nu_4$ & $30.3468$ & $30.3468$ & $30.3465$ \\ \hline
            $\nu_5$       & $46.9448$ & $46.9444$ & $46.9283$ \\ 
            $\tilde\nu_5$ & $46.9497$ & $46.9492$ & $46.9300$ \\ \hline
        \end{tabular}
    \end{subtable}%
    \begin{subtable}{0.5\linewidth}
        \caption{Outputs $y_2(t)$ and $\tilde y_2(t)$.}\label{tab:freq.op2} \centering
        \begin{tabular}{|c||c|c|c|} \hline
                          & $d=0.025$ & $d=1$     & $d=10$    \\ \hline
            $\nu_1$       & $4.3185$  & $4.3186$  & $4.3223$  \\ 
            $\tilde\nu_1$ & $4.3186$  & $4.3186$  & $4.3223$  \\ \hline
            $\nu_2$       & $12.6358$ & $12.6445$ & $13.5575$ \\ 
            $\tilde\nu_2$ & $12.6360$ & $12.6447$ & $13.5575$ \\ \hline
            $\nu_3$       & $18.5019$ & $18.4936$ & $17.6862$ \\ 
            $\tilde\nu_3$ & $18.5025$ & $18.4941$ & $17.6880$ \\ \hline
            $\nu_4$       & $30.3447$ & $30.3445$ & $30.3272$ \\ 
            $\tilde\nu_4$ & $30.3468$ & $30.3466$ & $30.3290$ \\ \hline
            $\nu_5$       & $46.9448$ & $46.9431$ & $46.7630$ \\ 
            $\tilde\nu_5$ & $46.9496$ & $46.9482$ & $46.7701$ \\ \hline
        \end{tabular}
    \end{subtable} \vspace{-2ex}
\end{table}

\section{Conclusion}
We conclude from the numerical investigations that the Timoshenko and Euler--Bernoulli models with an attached mass-spring actuator are in good agreement in terms of the localization of transfer functions poles at relatively low frequencies.
As observed in~\cite{RB2007,MM2005,ZY2020}, the effects of shear deformation of beam elements on vibration behavior become particularly significant in the medium- and high-frequency ranges, even when the beam’s length-to-width ratio is sufficiently large.
Therefore, further asymptotic analysis of the eigenfrequency distribution is required.

Changes in the damping coefficient slightly influence the frequency distribution, as confirmed by the data in Table~\ref{tab:freq}, but significantly affect the magnitude of the transfer function close to the modal frequencies, as illustrated in Figures~\ref{fig:op1} and~\ref{fig:op2}.
Since accurate information about the damping characteristics is often not available in practical applications, we consider the problem of damping parameter identification to be a promising direction for future investigation.

In the present paper, we limited ourselves to considering the output measured at the point $\ell_0$, where the actuator is located.
It is also of interest to investigate hybrid systems with non-collocated actuator-sensor combinations.
Moreover, we studied the local dissipation effect, while the distributed damping cases~\cite{ADW1986,LZ2018} are of interest for future research, particularly in the context of reduced-order modeling and subsequent control design.

\paragraph{Acknowledgment}
The first author gratefully acknowledges the funding by the European Regional Development Fund (ERDF) within the programme Research and Innovation~--- Grant Number ZS/2023/12/182138.
The second author was supported by a grant from the Simons Foundation (PD-Ukraine-00010584, Julia Kalosha).


\begin{thebibliography}{10}

\bibitem{AR2021}
Ahmed, A.M., Rifai, A.M.: \uppercase{E}uler--\uppercase{B}ernoulli and
  \uppercase{T}imoshenko beam theories analytical and numerical comprehensive
  revision.
\newblock European Journal of Engineering and Technology Research
  \textbf{6}(7), 20--32 (2021).
\newblock \url{doi: 10.24018/ejeng.2021.6.7.2626}

\bibitem{E2019}
Elishakoff, I.: Exact solution of \uppercase{T}imoshenko--\uppercase{E}hrenfest
  equations.
\newblock In: Handbook on \uppercase{T}imoshenko--\uppercase{E}hrenfest Beam
  and \uppercase{U}flyand--\uppercase{M}indlin Plate Theories, chap.~2, pp.
  107--138. World Scientific Book (2019).
\newblock \url{doi: 10.1142/9789813236523_0002}

\bibitem{T1921}
Timoshenko, S.P.: On the correction for shear of the
  differential equation for transverse vibrations of prismatic bars.
\newblock The London, Edinburgh, and Dublin Philosophical Magazine and Journal
  of Science \textbf{41}, 744--746 (1921).
\newblock \url{doi: 10.1080/14786442108636264}

\bibitem{T1974}
Timoshenko, S.P., Young, D.H., Weaver, W.: Vibration Problems in Engineering.
\newblock 4th Edition: Wiley, Chichester (1974).
\newblock \url{doi: 10.1017/S0001924000035041}

\bibitem{KS2000}
Krabs, W., Sklyar, G.M.: On the Stabilizability of a Slowly Rotating \uppercase{T}imoshenko Beam. Z.~Anal. Anwend. \textbf{19}(1), 131--145 (2000). \url{doi: 10.4171/ZAA/943}

\bibitem{LBMF1996}
Lueschen, G.G.G., Bergman, L.A., Mc\uppercase{F}arland, D.M.: Green's functions
  for uniform \uppercase{T}imoshenko beams.
\newblock Journal of Sound and Vibration \textbf{194}, 93--102 (1996).
\newblock \url{doi: 10.1006/jsvi.1996.0346}

\bibitem{TNC1953}
Traill-Nash, R.W., Collar, A.R.: The effects of shear flexibility and rotatory
  inertia on the bending vibrations of beams.
\newblock The Quarterly Journal of Mechanics and Applied Mathematics
  \textbf{6}, 186--222 (1953).
\newblock \url{doi: 10.1093/qjmam/6.2.186}

\bibitem{AWS1997}
Aldraihem, O.J., Wetherhold, R.C., Singh, T.: Distributed control of laminated
  beams: \uppercase{T}imoshenko theory vs.
  \uppercase{E}uler--\uppercase{B}ernoulli theory.
\newblock Journal of Intelligent Material Systems and Structures \textbf{8}(2),
  149--157 (1997).
\newblock \url{doi: 10.1177/1045389X9700800205}

\bibitem{K2021}
Ko\c{c}, M.A.: Finite element and numerical vibration analysis of a
  \uppercase{T}imoshenko and \uppercase{E}uler--\uppercase{B}ernoulli beams
  traversed by a moving high-speed train.
\newblock Journal of the Brazilian Society of Mechanical Sciences and
  Engineering \textbf{43}(165) (2021).
\newblock \url{doi: 10.1007/s40430-021-02835-7}

\bibitem{ZTS2020}
Zhang, X., Thompson, D., Sheng, X.: Differences between
  \uppercase{E}uler--\uppercase{B}ernoulli and \uppercase{T}imoshenko beam
  formulations for calculating the effects of moving loads on a periodically
  supported beam.
\newblock Journal of Sound and Vibration \textbf{481}, 115,432 (2020).
\newblock \url{doi: 10.1016/j.jsv.2020.115432}

\bibitem{P2005}
Park, J.: Transfer function methods to measure dynamic mechanical properties of
  complex structures.
\newblock Journal of Sound and Vibration \textbf{288}, 57--79 (2005).
\newblock \url{doi: 10.1016/j.jsv.2004.12.019}

\bibitem{RB2007}
Ruge, P., Birk, C.: A comparison of infinite \uppercase{T}imoshenko and
  \uppercase{E}uler--\uppercase{B}ernoulli beam models on \uppercase{W}inkler
  foundation in the frequency- and time-domain.
\newblock Journal of Sound and Vibration \textbf{304}, 932--947 (2007).
\newblock \url{doi: 10.1016/j.jsv.2007.04.001}

\bibitem{W1995}
Wang, C.M.: \uppercase{T}imoshenko beam-bending solutions in terms of
  \uppercase{E}uler--\uppercase{B}ernoulli solutions.
\newblock Journal of Engineering Mechanics \textbf{121}, 763--765 (1995).
\newblock \url{doi: 10.1061/(ASCE)0733-9399(1995)121:6(763)}

\bibitem{KZ2022}
Kalosha, J., Zuyev, A.: Asymptotic stabilization of a flexible beam with
  attached mass.
\newblock Ukrainian Mathematical Journal \textbf{73}, 1537--1550 (2022).
\newblock \url{doi: 10.1007/s11253-022-02012-6}

\bibitem{KZB2021}
Kalosha, J., Zuyev, A., Benner, P.: On the eigenvalue distribution for a beam
  with attached masses.
\newblock In: G.~Sklyar, A.~Zuyev (eds.) Stabilization of Distributed Parameter
  Systems: Design Methods and Applications, \emph{SEMA SIMAI Springer Series},
  vol.~2, pp. 43--56. Springer International Publishing, Cham (2021).
\newblock \url{doi: 10.1007/978-3-030-61742-4_3}

\bibitem{ZS05}
Zuyev, A., Sawodny, O.: Stabilization of a flexible manipulator model with passive joints.
\newblock IFAC Proceedings Volumes (38), 784--789 (2005).
\newblock \url{doi: 10.3182/20050703-6-CZ-1902.00531}

\bibitem{ZS2007}
Zuyev, A., Sawodny, O.: Stabilization and observability of a rotating
  \uppercase{T}imoshenko beam model.
\newblock Mathematical Problems in Engineering (1), 1--19 (2007).
\newblock \url{doi: 10.1155/2007/57238}

\bibitem{B2013}
Bartecki, K.: A general transfer function representation for a class of
  hyperbolic distributed parameter systems.
\newblock International Journal of Applied Mathematics and Computer Science
  \textbf{23}(2), 291--307 (2013).
\newblock \url{doi: 10.2478/amcs-2013-0022}

\bibitem{JB2006}
Jovanovic, M., Bamieh, B.: A formula for frequency responses of distributed
  systems with one spatial variable.
\newblock Systems Control Letters \textbf{55}(1), 27--37 (2006).
\newblock \url{doi: 10.1016/j.sysconle.2005.04.014}

\bibitem{GLO2017}
Guiver, C., Logemann, H., Opmeer, M.R.: Transfer functions of
  infinite-dimensional systems: positive realness and stabilization.
\newblock Mathematics of Control, Signals, and Systems \textbf{29}(20), 1--61
  (2017).
\newblock \url{doi: 10.1007/s00498-017-0203-z}

\bibitem{M2024}
Malti, R., Rapai{\'c}, M.R., Turkulov, V.: A unified framework for exponential
  stability analysis of irrational transfer functions in the parametric space.
\newblock Annual Reviews in Control \textbf{57}, 100,935 (2024).
\newblock \url{doi: 10.1016/j.arcontrol.2024.100935}

\bibitem{WF2015}
Wo{\'z}niak, J., Firkowski, M.: Note on the Stability of a Slowly Rotating \uppercase{T}imoshenko Beam with Damping.
\newblock Advances in Applied Mathematics and Mechanics. \textbf{7}(6), 736--753 (2015).
\newblock \url{doi: doi:10.4208/aamm.2014.m634}

%\bibitem{BBF2014}
%Baur, U., Benner, P., Feng, L.: Model order reduction for linear and nonlinear
%  systems: A system-theoretic perspective.
%\newblock Archives of Computational Methods in Engineering \textbf{21},
%  331--358 (2014).
%\newblock \url{doi: 10.1007/s11831-014-9111-2}
%
%\bibitem{DO2014}
%Demir, O., \uppercase{{\"O}}zbay, H.: On reduced order modeling of flexible
%  structures from frequency response data.
%\newblock In: 2014 European Control Conference (ECC), pp. 1133--1138. IEEE,
%  Strasbourg, France (2014).
%\newblock \url{doi: 10.1109/ECC.2014.6862456}
%
%\bibitem{KGA2021}
%Karachalios, D.S., Gosea, I.V., Antoulas, A.C.: The \uppercase{L}oewner
%  framework for system identification and reduction.
%\newblock In: P.~Benner, S.~Grivet-Talocia, A.~Quarteroni, G.~Rozza,
%  W.~Schilders, L.M. Silveira (eds.) Model Order Reduction: Volume 1: System-
%  and Data-Driven Methods and Algorithms, pp. 181--228. de Gruyter (2021).
%\newblock \url{doi: 10.1515/9783110498967-006}

\bibitem{LHW1999}
Lew, J., Hyde, C., Wade, M.: Damage detection of flexible beams using transfer
  function correlation.
\newblock In: Proceedings of the ASME 1999 Design Engineering Technical
  Conferences, \emph{International Design Engineering Technical Conferences and
  Computers and Information in Engineering Conference}, vol. Volume 7A: 17th
  Biennial Conference on Mechanical Vibration and Noise, pp. 905--912. Las
  Vegas, Nevada, USA (1999).
\newblock \url{doi: 10.1115/DETC99/VIB-8373}

\bibitem{CM2009}
Curtain, R., Morris, K.: Transfer functions of distributed parameter systems: A
  tutorial.
\newblock Automatica \textbf{45}, 1101--1116 (2009).
\newblock \url{doi: 10.1016/j.automatica.2009.01.008}

\bibitem{LA1993}
Lenz, K., \uppercase{{\"O}}zbay, H.: Analysis and robust control techniques for
  an ideal flexible beam.
\newblock In: C.~Leondes (ed.) Multidisciplinary Engineering Systems: Design
  and Optimization Techniques and their Application, \emph{Control and Dynamic
  Systems}, vol.~57, pp. 369--421. Academic Press (1993).
\newblock \url{doi: 10.1016/B978-0-12-012757-3.50015-4}

\bibitem{ZG2023}
Zuyev, A., Gosea, I.V.: Approximating a flexible beam model in the
  \uppercase{L}oewner framework.
\newblock In: 2023 European Control Conference (ECC), pp. 1--7. IEEE,
  Bucharest, Romania (2023).
\newblock \url{doi: 10.23919/ECC57647.2023.10178203}

\bibitem{DS2013}
Dym, C.L., Shames, I.H.: Solid Mechanics: A Variational Approach, Augmented
  Edition.
\newblock Springer New York, NY (2013).
\newblock \url{doi: 10.1007/978-1-4614-6034-3}

\bibitem{GT1991}
Gere, J.M., Timoshenko, S.P.: Mechanics of materials. 3rd Ed.
\newblock Chapman and Hall. Ltd., London. England (1991)

\bibitem{LCR2017}
Lenci, S., Clementi, F., Rega, G.: Comparing nonlinear free vibrations of
  \uppercase{T}imoshenko beams with mechanical or geometric curvature definition.
\newblock Procedia IUTAM \textbf{20}, 34--41 (2017).
\newblock \url{doi: 10.1016/j.piutam.2017.03.006}

\bibitem{MM2005}
Mei, C., Mace, B.R.: Wave reflection and transmission in \uppercase{T}imoshenko
  beams and wave analysis of \uppercase{T}imoshenko beam structures.
\newblock Journal of Vibration and Acoustics \textbf{127}(4), 382--394 (2005).
\newblock \url{doi: 10.1115/1.1924647}

\bibitem{ZY2020}
Zhang, Y., Yang, B.: Medium-frequency vibration analysis of
  \uppercase{T}imoshenko beam structures.
\newblock International Journal of Structural Stability and Dynamics
  \textbf{20}(13), 2041,009 (2020).
\newblock \url{doi: 10.1142/S0219455420410096}

\bibitem{ADW1986}
Alberts, T.E., Dickerson, S.L., Book, W.J.: On the transfer function modeling
  of flexible structures with distributed damping.
\newblock Robotics: Theory and Applications, ASME pp. 23--30 (1986)

\bibitem{LZ2018}
Liu, Z., Zhang, Q.: Stability and regularity of solution to the
  \uppercase{T}imoshenko beam equation with local
  \uppercase{K}elvin--\uppercase{V}oigt damping.
\newblock SIAM Journal on Control and Optimization \textbf{56}, 3919--3947
  (2018).
\newblock \url{doi: 10.1137/17M1146506} 

\end{thebibliography}
\end{document}